\documentclass[letterpaper, 10 pt, conference]{IEEEconf}  

\IEEEoverridecommandlockouts                              %
\usepackage{cite}
\usepackage{url}
\usepackage{amsmath,amssymb,amsfonts}
\usepackage{multirow}
\usepackage{algorithmic}
\usepackage{graphicx}
\usepackage{textcomp}
\usepackage{epsfig}
\usepackage{algorithmic}
\usepackage{epstopdf}
\usepackage[ruled]{algorithm}
\usepackage{enumerate}
\usepackage{float}
\usepackage{graphics,graphicx}
\usepackage{subfigure}
\usepackage{booktabs}  

\usepackage{array}
\usepackage{mathtools}
\usepackage[justification=centering]{caption}
\allowdisplaybreaks[4]
\usepackage{url}
\DeclareGraphicsExtensions{.png}
\graphicspath{{./Pics/}}
\usepackage{color}
\usepackage{xspace}

\newtheorem{lemma}{Lemma}
\newtheorem{theorem}{Theorem}

\newtheorem{assumption}{Assumption}
\newtheorem{proposition}{Proposition}

\title{\LARGE \bf A Zeroth-Order Proximal Algorithm for Consensus Optimization
}

\author{Chengan Wang, Zichong Ou and
	Jie Lu
\thanks{C. Wang, Z. Ou and J. Lu are with the School of Information Science and Technology, ShanghaiTech University, 201210, Shanghai, China. Email:{\tt \{wangcha, ouzch, lujie\}@shanghaitech.edu.cn}.}\\
}

\begin{document}

\maketitle
\thispagestyle{empty}
\pagestyle{empty}

\begin{abstract} 
This paper considers a consensus optimization problem, where all the nodes in a network, with access to the zeroth-order information of its local objective function only, attempt to cooperatively achieve a common minimizer of the sum of their local objectives. To address this problem, we develop ZoPro, a zeroth-order proximal algorithm, which incorporates a zeroth-order oracle for approximating Hessian and gradient into a recently proposed, high-performance distributed second-order proximal algorithm. We show that the proposed ZoPro algorithm, equipped with a dynamic stepsize, converges linearly to a neighborhood of the optimum in expectation, provided that each local objective function is strongly convex and smooth. Extensive simulations demonstrate that ZoPro converges faster than several state-of-the-art distributed zeroth-order algorithms and outperforms a few distributed second-order algorithms in terms of running time for reaching given accuracy. 
\end{abstract}


\section{INTRODUCTION}
This paper considers  a widely studied distributed optimization problem, i.e., consensus optimization, where all the nodes in a network aim at reaching a consensus that minimizes the sum of their local cost functions. This problem arises in many real-world applications such as distributed machine learning \cite{b1} and resource allocation \cite{b2}.

To date, a variety of distributed algorithms for convex consensus optimization have been proposed, in which each node only has access to certain information of its convex local cost function and can only communicate with its neighbors determined by the network topology. Most existing distributed optimization algorithms are first-order methods, which typically include primal (sub-)gradient methods \cite{b3,b4,b5,b6,b7} and dual/primal-dual (sub-)gradient methods \cite{b10,b11,b12}. These methods essentially require the nodes compute the (sub-)gradients of their primal or dual objectives. 

The second-order methods, such as the Decentralized Broyden-Fletcher-Goldfarb-Shanno method (D-BFGS) \cite{b16}, the Exact Second-Order Method (ESOM) \cite{b17}, the Decentralized Quadratically Approximated ADMM (DQM) \cite{b18} and the Second-Order Proximal Algorithm (SoPro) \cite{b20}, employ the objective Hessian matrices in addition to the objective gradients, potentially leading to faster convergence due to their more accurate approximations of some global objectives. Recently, Wu \textit{et al}. have proposed a distributed second-order proximal algorithm called SoPro \cite{b20}. SoPro originates from the classic Method of Multipliers \cite{b33}, while it replaces the augmented Lagrangian function of the problem with its second-order approximation and introduces a separable quadratic proximal term to decouple the problem. SoPro achieves linear convergence under strong convexity and smoothness, and exhibits superior practical performance.

All of the aforementioned algorithms are required to compute at least the first-order information (i.e., (sub-)gradients) and even the second-order information (i.e., Hessian matrices) of the problem. However, these pieces of information could be unavailable or too time-consuming to obtain in big data and large-scale network scenarios. Under such circumstances, zeroth-order algorithms \cite{b22,b23,b24,b25,b26,b27,b31,b32} are effective approaches, whose updates only involve sampled function values instead of the exact (sub-)gradients and Hessian matrices. Various zeroth-order algorithms have been developed so far, which utilize zeroth-order estimators with different batch sizes to estimate the (sub-)gradients in some  distributed first-order methods, including zeroth-order gradient tracking method \cite{b22}, zeroth-order primal-dual methods \cite{b23,b24}, zeroth-order decentralized (sub-)gradient descent \cite{b22}, zeroth-order method with approximate projection\cite{b26}, distributed zeroth-order projected gradient descent \cite{b27}, distributed randomized zeroth-order mirror descent method \cite{b32}, etc. To the best of our knowledge, zeroth-order Hessian estimators have barely been considered in distributed optimization. As distributed second-order algorithms often outperform the first-order ones in both accuracy and convergence rate, introducing zeroth-order oracles for Hessian estimation to distributed second-order methods is a promising direction.


In this paper, we propose a distributed \underline{z}eroth-\underline{o}rder \underline{pro}ximal algorithm, referred to as ZoPro, for solving convex consensus optimization. ZoPro replaces the exact objective gradients and Hessian matrices in the distributed second-order algorithm SoPro \cite{b20} with their zeroth-order estimates, so that it significantly reduces the computational cost of SoPro and is applicable to the scenarios where the objective gradients and Hessian matrices are inaccessible or too costly to compute. ZoPro also inherits the appealing convergence performance of SoPro. It is shown to achieve a linear rate of convergence to a neighborhood of the optimal solution when the objective functions are strongly convex and smooth. Finally, the numerical experiments demonstrate that ZoPro outperforms a few state-of-the-art distributed zeroth-order methods in convergence speed and enjoys shorter running time than several well-noted second-order methods.

The rest of the paper is organized as follows: Section II describes the problem formulation, Section III develops the proposed ZoPro algorithm, Section IV provides convergence result, Section V presents the simulation results, and Section VI concludes the paper.

\textbf{Notations and definitions:} For any differentiable function $f:\mathbb{R}^d\rightarrow\mathbb{R}$, $\nabla f(x)$ represents its gradient at $x\in\mathbb{R}^d$ and if $f$ is twice-differentiable, $\nabla^2 f(x)$ represents its Hessian matrix. $\mathbf{O}_d$ and $\mathbf{I}_d$ represent the $d\times d$ zero matrix and identity matrix, respectively, and $\mathbf{0}_d$ and $\mathbf{1}_d$ represent the $d$-dimensional all-zero and all-one vectors, respectively. Define $\mathbb{Z}^+$ and $\mathbb{N}$ as the sets of positive integers and non-negative integers, respectively. Also, $\otimes$ denotes the Kronecker product, $\|\cdot\|$ denotes the $\mathcal{L}_2$ norm and $\left\langle\cdot,\cdot\right\rangle$ denote the inner product. Besides, $\text{diag}(A_1,\dots,A_n)$ represents the block diagonal matrix consisting of the diagonal blocks $A_1,\dots,A_n$. $[P]_{ij}$ denotes the $(i,j)$-entry of matrix $P$. 
Given $A=A^T\in\mathbb{R}^{d\times d}$ and $\mathbf{x}\in\mathbb{R}^d$, $\|\mathbf{x}\|^2_A=\mathbf{x}^T A\mathbf{x}$. $\lambda_{\min}(A)$ and $\lambda_{\max}(A)$ represent $A$'s smallest and largest eigenvalues, respectively. $A^{\dagger}$ denotes $A$'s pseudo-inverse and $A^\perp$ represents the orthogonal complement of $A$. $v\sim\mathcal{N}(\mu,\Sigma)$ represents a Gaussian random vector $v$ with mean $\mu$ and covariance matrix $\Sigma$. A function $f:\mathbb{R}^d\rightarrow\mathbb{R}$ is $\mu$-strongly convex if $f$ is differentiable and 
\vspace{-0.15cm}\begin{equation*}
	f(y)\geq f(x)+\nabla f(x)^T(y-x)+\frac{\mu}{2}\|y-x\|^2\ \forall x,y\in\mathbb{R}^d
	\vspace{-0.1cm}\end{equation*}
for some $\mu>0$. $f$ is $L$-smooth if $f$ is differentiable and
\vspace{-0.15cm}\begin{equation*}
	\|\nabla f(x)-\nabla f(y)\|\leq L\|x-y\|,\ \forall x,y\in\mathbb{R}^d
\vspace{-0.15cm}\end{equation*}
for some $L>0$. Finally, the directional derivative of $f$ at point $x$ along direction $d$ is denoted by $f'(x;d)=\lim_{\alpha\rightarrow 0}\frac{f(x+\alpha d)-f(x)}{\alpha}$.

\section{PROBLEM FORMULATION}\label{section problem form}
\label{sec:preliminaries}
We consider solving
\begin{equation}\label{objectiveproblem}
	\underset{x\in\mathbb{R}^{d}}{\operatorname{minimize}}\sum_{i\in\mathcal{V}}f_i(x)
\end{equation}
over a network modeled as a connected and undirected graph $\mathcal{G}=(\mathcal{V},\mathcal{E})$, where $\mathcal{V}=\left\{1,2,\dots,N\right\}$ is the set of nodes and $\mathcal{E}\subseteq\left\{\left\{i,j\right\}\subseteq\mathcal{V}\times\mathcal{V}\ |\ i\neq j\right\}$ is the set of the bidirectional links. For each $i\in\mathcal{V}$, we denote the set of its neighbors by $\mathcal{N}_i=\left\{j\in\mathcal{V}\ |\ \{i,j\}\in\mathcal{E}\right\}$. Here,  $f_i:\mathbb{R}^d\rightarrow\mathbb{R}$ is the local cost/objective function associated with node $i\in\mathcal{V}$ and each node $i$ only communicates with its neighbors in $\mathcal{N}_i$.

Let $x_i\in\mathbb{R}^d$ be node $i$’s local copy of the global optimization variable $x$ and $\mathbf{x}$ be the concatenation of all the $x_i$'s, i.e. $\mathbf{x}=(x_1^T,\dots,x_N^T)^T\in\mathbb{R}^{Nd}$. Let $P=P^T$ be a weight matrix corresponding to the network $\mathcal{G}$ given by 
\begin{align*}
	[P]_{i j}=\left\{\begin{array}{ll}
		\sum_{s \in \mathcal{N}_{i}} p_{i s}, & i=j, \\
		-p_{i j}, & j \in \mathcal{N}_{i}, \\
		0, & \text {otherwise,}
	\end{array} \quad \forall i, j \in \mathcal{V},\right.
\end{align*}
where $p_{ij}=p_{ji}>0\ \forall\left\{i,j\right\}\in\mathcal{E}$.

Due to the fact that $\mathcal{G}$ is connected, the null space of $P$ is span$\left\{\mathbf{1}_N\right\}$, so that we can rewrite problem (\ref{objectiveproblem}) as
\vspace{-0.1cm}\begin{align}
	&\underset{x\in\mathbb{R}^{Nd}}{\operatorname{minimize}}\ f(\mathbf{x})=\sum_{i\in\mathcal{V}}f_i(x_i) \nonumber\\ 
	&\text{subject to } W^{\frac{1}{2}}\mathbf{x}=\mathbf{0}_{Nd}, \label{consensus_problem}
\end{align}
where $W=P\otimes\mathbf{I}_d\succeq  \mathbf{O}_{Nd}$ and the equality constraint means that $x_1,\dots,x_n$ are identical \cite{b20}.

We impose the following assumption on problem (\ref{objectiveproblem}).
\begin{assumption}\label{assump1}
	Each $f_i$ is $m_i$-strongly convex, twice continuously differentiable and $M_i$-smooth, where $m_i$, $M_i>0$.
\end{assumption}

Note that Assumption \ref{assump1} guarantees the uniqueness of the optimal solution $x^*$ to problem (\ref{objectiveproblem}).

\section{ALGORITHM DEVELOPMENT}\label{section algorithm design}

In this section, we develop a distributed zeroth-order algorithm for solving problem (\ref{consensus_problem}). 

\subsection{SoPro Algorithm}
We first quickly review the second-order proximal (SoPro) algorithm proposed in \cite{b20}.

SoPro solves (\ref{consensus_problem}) in a primal-dual fashion as follows: 
\begin{align}
	\label{primal_iteration2}
	\mathbf{x}^{k+1}&=\mathbf{x}^{k}-(\triangledown^2f(\mathbf{x}^k)+D)^{-1}(\triangledown f(\mathbf{x}^k)+\rho W\mathbf{x}^k+\mathbf{q}^k),\displaybreak[0]\\
	\label{dual_iteration2}
	\mathbf{q}^{k+1}&=\mathbf{q}^k+\rho W\mathbf{x}^{k+1},
\end{align}
with the initialization $\mathbf{q}^0=\mathbf{0}_{Nd}$, where $\mathbf{x}^k$ is the global primal variable and $\mathbf{q}^k=W^{\frac{1}{2}}\mathbf{v}^k$ is a change of variable with $\mathbf{v}^k$ being the dual variable associated with the constraint in (\ref{consensus_problem}). The primal update (\ref{primal_iteration2}) intends to minimize an augmented-Lagrangian-like function constructed in the following way: The augmented Lagrangian function of (\ref{consensus_problem}) with penalty $\frac{\rho}{2}\|W^{\frac{1}{2}}\mathbf{x}\|^2$, $\rho>0$ is first replaced by its second-order approximation at $\mathbf{x}^k$ to reduce computational cost. Then, a separable quadratic proximal term $\frac{1}{2}(\mathbf{x}-\mathbf{x}^k)^T(\nabla^2f(\mathbf{x}^k)+D)(\mathbf{x}-\mathbf{x}^k)$ is employed as a substitute for the non-separable term in the above approximate augmented Lagrangian function to enable fully distributed implementation, where $D=\text{diag}(D_1,\dots,D_N)$ is a symmetric block diagonal matrix with each $D_i\in\mathbb{R}^{d\times d}$ satisfying $\nabla^2 f(\mathbf{x}^k)+D\succ\mathbf{O}_{Nd}$. The dual update (\ref{dual_iteration2}) emulates dual gradient ascent.

Observe that the updates (\ref{primal_iteration2}) and (\ref{dual_iteration2}) of SoPro require calculating accurate first-order and second-order information of the objective function $f$, which could be a tough challenge when handling big data and large-scale problems. The high computational complexity of SoPro motivates the development of zeroth-order oracles for efficiently estimating the gradients and Hessian in SoPro.

\subsection{Zeroth-order Oracle}
Next, we provide a zeroth-order oracle for estimating the gradients and Hessian matrices in SoPro's updates (\ref{primal_iteration2}) and (\ref{dual_iteration2}). To do so, consider the following smoothed approximation of the objective function $f$:
\begin{equation}
	f_\mu(\mathbf{x})\triangleq\frac{1}{(2\pi)^{Nd/2}}\int_{\mathbb{R}^{Nd}}f(\mathbf{x}+\mu\mathbf{u})e^{-\frac{\|\mathbf{u}\|^2}{2}}d\mathbf{u},
\end{equation}
where $\mathbf{u}\sim(0,\mathbf{I}_{Nd})\in\mathbb{R}^{Nd}$ is a Gaussian random vector and $\mu>0$ is a parameter to control the smoothness level \cite[Section 2]{b34}. Note that the smoothed approximation $f_\mu$ is guaranteed to be differentiable. We will show more properties of $f_\mu$ in Section \ref{convergence analysis}. Let $\widetilde{g}_{\mu}$ and $\widetilde{H}_{\mu}$ represent the gradient and Hessian estimation of the above Gaussian smoothing function $f_\mu$, which are defined according to \cite{b31} as follows:
\vspace{-0.15cm}\begin{equation}\label{gradient_estimate}
	\widetilde{g}_\mu(\mathbf{x})=\frac{1}{b}\sum_{j=1}^b\frac{f(\mathbf{x}+\mu u_j)-f(\mathbf{x})}{\mu}u_j.
\end{equation}
\begin{equation}  
	\label{hessian_estimate_i}
	\widetilde{H}_{\mu}(\mathbf{x}) =\text{diag}\left(\widetilde{H}_{\mu,1}(x_1),\dots,\widetilde{H}_{\mu,N}(x_N)\right),
	\vspace{-0.15cm}\end{equation}
\vspace{-0.2cm}\begin{align*}
	&\widetilde{H}_{\mu,i}(x_i) \nonumber\\
	=&\frac{1}{b}\sum_{j=1}^b\frac{f_i(x_i+\mu u_j)+f_i(x_i-\mu u_j)-2f_i(x_i)}{2\mu^2}u_j u_j^T,\nonumber\\
	&\forall i=1,\dots,N,
	\vspace{-0.2cm}\end{align*}
where $b\in\mathbb{Z}^+$ is the batch size, $u_j\sim\mathcal{N}(0,\mathbf{I}_{d})\in\mathbb{R}^{d}$, $j=1,\dots,b$ are Gaussian random vectors. The zeroth-order oracle (\ref{gradient_estimate})--(\ref{hessian_estimate_i}) for gradient and Hessian matrix estimation only needs to sample $(2b+1)N$ points from the local objective functions, which is much less costly than computing the exact gradient and Hessian matrix. It can be verified that the zeroth-order gradient estimation (\ref{gradient_estimate}) is an unbiased estimator for $\nabla f_\mu$ \cite{b25}, i.e., $\mathbf{E}\left[\widetilde{g}_\mu(\mathbf{x})\right]=\nabla f_\mu(\mathbf{x})$.

\subsection{ZoPro Algorithm}
In this subsection, we incorporate the zeroth-order oracle (\ref{gradient_estimate})--(\ref{hessian_estimate_i}) into SoPro, yielding a zeroth-order proximal algorithm, referred to as ZoPro.

We first replace $\nabla f(\mathbf{x}^k)$ and $\nabla^2 f(\mathbf{x}^k)$ in (\ref{primal_iteration2}) and (\ref{dual_iteration2}) with $\widetilde{g}_\mu(\mathbf{x}^k)$ and $\widetilde{H}_\mu(\mathbf{x}^k)$, which gives
\vspace{-0.15cm}\begin{align}
	\label{zeroth_primal_iteration}
	\mathbf{x}^{k+1}&=\mathbf{x}^k-(\widetilde{H}_{\mu}(\mathbf{x}^k)+D)^{-1}(\widetilde{g}_{\mu}(\mathbf{x}^k)+\rho W\mathbf{x}^k+\mathbf{q}^k),\\
	\label{zeroth_dual_iteration}
	\mathbf{q}^{k+1}&=\mathbf{q}^k+\rho W\mathbf{x}^{k+1},
	\vspace{-0.1cm}\end{align}
where, similar to SoPro, $D=\text{diag}(D_1,\dots,D_N)\in\mathbb{R}^{Nd\times Nd}$ is a symmetric block diagonal matrix such that $\widetilde{H}_{\mu}(\mathbf{x})+D\succ\mathbf{O}_{Nd}\ \forall \mathbf{x}\in\mathbb{R}^{Nd}$, or equivalently, $\widetilde{H}_{\mu,i}(x_i)+D_i\succ\mathbf{0}_{d}$ $\forall x\in\mathbb{R}^d$ $\forall i=1,\dots,N$. The starting point $\mathbf{q}^0$ is set to $\mathbf{q}^0\in S^\perp$, where $S=\left\{\mathbf{x}\in\mathbb{R}^{Nd}|x_1=\dots=x_N\right\}$ and $S^\perp=\left\{\mathbf{x}\in\mathbb{R}^{Nd}|x_1+\dots+x_N=\mathbf{0}_d\right\}$, so that $\mathbf{q}^k\in S^\perp$ $\forall k\geq 0$ due to (\ref{zeroth_dual_iteration}). For simplicity, we set $\mathbf{q}^0=\mathbf{0}_{Nd}$.

Since $\widetilde{g}_\mu(\mathbf{x}^k)$ and $\widetilde{H}_\mu(\mathbf{x}^k)$ are only estimated values of $\nabla f(\mathbf{x}^k)$ and $\nabla^2 f(\mathbf{x}^k)$, (\ref{zeroth_primal_iteration}) and (\ref{zeroth_dual_iteration}) may not converge to the exact optimum like SoPro. To overcome this issue, we introduce a backtracking line search strategy with a dynamic stepsize to bound the sequence $\left\{\mathbf{x}^k\right\}$. We set the search direction to $\mathbf{d}^k=-(\widetilde{H}_{\mu}(\mathbf{x}^k)+D)^{-1}(\widetilde{g}_{\mu}(\mathbf{x}^k)+\rho W\mathbf{x}^k+\mathbf{q}^k)$, and then modify (\ref{zeroth_primal_iteration}) to 
\vspace{-0.2cm}\begin{align}
	\label{zeroth_primal_backtracking_iteration}
	\mathbf{x}^{k+1}&=\mathbf{x}^k+A^k\mathbf{d}^k.
	\vspace{-0.2cm}\end{align}
Here, $A^k=\text{diag}(\alpha_1^k,\alpha_2^k,\dots,\alpha_N^k)\otimes\mathbf{I}_d\in\mathbb{R}^{Nd\times Nd}$ and $\alpha_i^k,\ i\in\mathcal{V}$ is the local stepsize of node $i$ determined by the Armijo condition \cite[Eq. (1)]{b35}, i.e. $f_i(x_i^k+\alpha_i^kd_i^k)\leq f(x_i^k)+c\alpha_i^kf_i'(x_i^k;d_i^k)$, where $c\in(0,1)$ is the stepsize control parameter, $f_i'(x_i^k;d_i^k)$ is the directional derivative of $f_i$ at $x_i^k$ along node $i$'s local search direction $d_i^k$ and $d_i^k$ is the $i$-th d-dimensional block of $\mathbf{d}^k$.

The primal update (\ref{zeroth_primal_backtracking_iteration}) and the dual update (\ref{zeroth_dual_iteration}) with initialization $\mathbf{q}^0=\mathbf{0}_{Nd}$ constitute a zeroth-order proximal algorithm, referred to as ZoPro, whose distributed implementation is described in Algorithm 1.

\begin{algorithm}
	\caption{Zeroth-Order Proximal Algorithm (ZoPro)}
	\begin{algorithmic}[1]
		\STATE \textbf{Initialization:}\\
		All the nodes agree on the batch size $b\in\mathbb{Z}^+$, the smoothness parameter $\mu>0$, the penalty parameter $\rho>0$, and the stepsize control parameter $c\in(0,1)$. Generate $b$ random vectors $u_j\sim\mathcal{N}(0,\mathbf{I}_d)$ $\forall j=1,\dots,b$.
		\STATE Each node $i\in\mathcal{V}$ chooses $D_i$ such that $\widetilde{H}_{\mu,i}(x)+D_i\succ\mathbf{O}_{d}\ \forall x\in\mathbb{R}^{d}$ and sets the initial stepsize $\alpha_i^0=1$.
		\STATE Every pair of neighboring nodes $\{i,j\}\in\mathcal{E}$ set $p_{ij}=p_{ji}$ to some positive value.
		\STATE Each node $i\in\mathcal{V}$ sets $x_{i}^0\in\mathbb{R}^d$ arbitrarily and $q_i^0=\mathbf{0}_d$. Then, it sends $x_{i}^0$ to every neighbor $j\in\mathcal{N}_i$.
		\STATE Upon receiving $x_j^0$ $\forall j\in\mathcal{N}_i$, each node $i\in\mathcal{V}$ sets $y_i^0=\sum_{j\in\mathcal{N}_i}p_{ij}(x_i^0-x_j^0)$.
		\FOR{$k\geq 0$}
		\STATE Each node $i\in\mathcal{V}$ computes Hessian estimate $\widetilde{H}_{\mu,i}(x_i^k)\linebreak[4]=b^{-1}\sum_{j=1}^b\frac{f_i(x_i^k+\mu u_j)+f_i(x_i^k-\mu u_j)-2f_i(x_i^k)}{2\mu^2}u_ju_j^T$.
		\STATE Each node $i\in\mathcal{V}$ computes gradient estimate $\widetilde{g}_{\mu,i}(x_i^k)=\frac{1}{b}\sum_{j=1}^b\frac{f_i(x_i^k+\mu u_j)-f_i(x_i^k)}\mu u_j$.
		\STATE Each node $i\in\mathcal{V}$ computes the search direction $d_i^k=-(\widetilde{H}_{\mu,i}(x_i^k)+D_i)^{-1}(\widetilde{g}_{\mu,i}(x_i^k)+\rho y_i^k +q_i^k)$.
		\STATE Each node $i\in\mathcal{V}$ determines the stepsize $\alpha^k_i$ such that $f_i(x_i^k+\alpha_i^kd_i^k)\leq f(x_i^k)+c\alpha_i^kf_i'(x_i^k;d_i^k)$.
		\STATE Each node $i\in\mathcal{V}$ updates $x_i^{k+1}=x_i^{k}+\alpha^k_id_i^k$ and sends $x_i^{k+1}$ to every neighbor $j\in\mathcal{N}_i$.
		\STATE Upon receiving $x_j^{k+1}$ $\forall j\in\mathcal{N}_i$, each node $i\in\mathcal{V}$ updates $y_i^{k+1}=\sum_{j\in\mathcal{N}_i}p_{ij}(x_i^{k+1}-x_j^{k+1})$ and $q_i^{k+1}=q_i^k+\rho y_i^{k+1}$.
		\ENDFOR
	\end{algorithmic}
\end{algorithm}

In Algorithm 1, each node $i$ maintains a local primal variable $x_i^k\in\mathbb{R}^d$ and a local dual variable $q_i^k\in\mathbb{R}^d$, which are the $i$-th d-dimensional block of $\mathbf{x}^k$ and $\mathbf{q}^k$. Also, we let it maintain an auxiliary variable $y_i^k\in\mathbb{R}^d$ such that $\mathbf{y}^k=\left((y_1^k)^T,\dots,(y_N^k)^T\right)^T=W\mathbf{x}^k$ for better presentation.

The existing zeroth-order distributed optimization methods such as distributed zeroth-order gradient tracking method \cite{b22}, distributed zeroth-order primal-dual method \cite{b23,b24}, distributed zeroth-order projected gradient descent \cite{b27} and distributed randomized zeroth-order mirror descent method \cite{b32} all use zeroth-order information to approximate the objective gradients only. In contrast, ZoPro includes zeroth-order estimates for both gradients and Hessian matrices. This may accelerate the convergence as ZoPro adopts potentially more accurate approximations of the global objective than other zeroth-order methods.

\section{CONVERGENCE ANALYSIS}\label{convergence analysis}
\vspace{-0.15cm}
This section provides the convergence analysis of ZoPro.

First, we analyze some properties of $f_i$ and $f_i'$.

\begin{proposition}\label{prop1}
	Let $f_i:\mathbb{R}^{d}\rightarrow\mathbb{R}$ be a L-smooth function and let $\left\{x_i^{k}\right\}$ be the sequence generated by $x_i^{k+1}=x_i^{k}+\alpha_i^{k}d_i^{k}$, where $\alpha_i^{k}$ is the stepsize determined by the backtracking line search and $d_i^k$ is the corresponding search direction. Denote the directional derivative of $f_i$ as $f_i'$. Then one of the following statements is
	true:\\
	(i) $f_i(x_i^k)\rightarrow -\infty$ as $k\rightarrow \infty$.\\
	(ii) The sequence $\left\{\|d_i^{k}\|\right\}$ diverges.\\
	(iii) For every infinite subsequence $J\subseteq \mathbb{N}$ for which $\left\{d_i^{k}: k\in J\right\}$ is bounded, we have
	\vspace{-0.2cm}\begin{equation*}
		\lim_{k\in J,k\rightarrow\infty} f_i'(x_i^k;d_i^k) = 0.
	\end{equation*}
\end{proposition}
\begin{proof}
See Appendix \ref{Proof of Proposition1}.
\end{proof}

From Proposition \ref{prop1}, we have $\lim_{k\rightarrow\infty}f_i'(x_i^k;d_i^k)=0$ for $f_i$, $i\in\mathcal{V}$, which indicates that either the gradient $\nabla f_i(x_i^k)$ and the search direction $d_i^k$ are orthogonal as $k\rightarrow\infty$ or $\nabla f_i(x_i^k)$ equals $0$ as $k\rightarrow\infty$, and both results can terminate the backtracking line search. From \cite[Theorem 3.2]{b36}, backtracking line search method guarantees convergence of the generated sequence $\left\{\mathbf{x}^k\right\}$. By continuous mapping theorem, $\left\{\nabla f(\mathbf{x}^k)\right\}$ is convergent, and thus $\left\{ \|\nabla f(\mathbf{x}^k)\|\right\}$ is bounded. For simplicity, we denote the upper bound as $K$.

Assumption \ref{assump1} implies $f$ is strongly convex for some $m\in\left(0,\min_{i}m_{i}\right]$ and smooth for some $M\geq\max_{i}M_{i}$. \cite[Eq. (9)]{b25} derives a bound of the difference between $\nabla f_\mu$ and $\widetilde{g}_\mu$.
\vspace{-0.15cm}\begin{align}
	&\mathbf{E}\left[\|\widetilde{g}_\mu(\mathbf{x})-\nabla f_\mu(\mathbf{x})\|^2\right]\leq\frac{2Nd\left(\mu^2M^2Nd+K^2\right)}{b}\triangleq G_1^2 \label{norm_bound}.
\vspace{-0.15cm}\end{align}  
Besides, the difference between $\nabla f_\mu$ and $\nabla f$ is also bounded in \cite[Eq. (28)]{b34} as follows:
\vspace{-0.15cm}\begin{align}
	\left\|\nabla f_\mu(\mathbf{x})-\nabla f(\mathbf{x})\right\|^2\leq \frac{\mu^2}{4}M^2(Nd+3)^3\triangleq G_2^2. \label{norm_bound_2}
\vspace{-0.15cm}\end{align}

Moreover, let $f_{\mu,i}$ be the smoothed approximation of function $f_i$, i.e. $f_{\mu,i}(x_i)=\frac{1}{(2\pi)^{d/2}}\int_{\mathbb{R}^d}f(\mathbf{x}+\mu u_i)e^{-\frac{\|u_i\|^2}{2}}du_i$, where, similarly,  $u_i\sim(0,\mathbf{I}_{d})\in\mathbb{R}^{d}$, $i=1,\dots,N$ is Gaussian random vectors. According to \cite[Section 2]{b34}, the smoothing function $f_{\mu,i}$ can  preserve all characteristics of $f_i$. For example, $f_{\mu,i}$ is guaranteed to be $m_i$-strongly convex and $M_i$-smooth if $f_i$ is $m_i$-strongly convex and $M_i$-smooth. $f_{\mu,i}$ is twice continuously differentiable if $f_i$ is twice continuously differentiable.
Assumption \ref{assump1} also implies $m_{i}\mathbf{I}_d\preceq\nabla^2 f_{i}(x)\preceq M_{i}\mathbf{I}_d$, $i=1,\dots,N$ $\forall x\in\mathbb{R}^d$. Let $\Lambda_m=\text{diag}(m_{1},m_{2},\dots,m_{N})\otimes \mathbf{I}_d\succ\mathbf{O}_{Nd}$ and  $\Lambda_M=\text{diag}(M_{1},M_{2},\dots,M_{N})\otimes \mathbf{I}_d\succ\mathbf{O}_{Nd}$. Besides, in order to provide a bound of the stepsize in backtracking line search, we define the smallest stepsize in the whole process as $\underline{\alpha}$, i.e. $\underline{\alpha}=\min_{i,k}\alpha_i^k>0$ for $i=1,\dots,N$ and $k=0,1,\dots$. We have $\underline{\alpha}\in(0,1]$ since $\alpha_i^0=1$. Let $R=\underline{\alpha}^{-1}(\frac{\Lambda_M+\Lambda _m}{2}+D)\in\mathbb{R}^{Nd\times Nd}$ and $Q=\text{diag}(\rho R,\mathbf{I}_{Nd})$. For simplicity, define $H^k=\widetilde{H}_\mu(\mathbf{x}^k)+D$. 

Similar to \cite[Eq. (3.1)]{b31}, we impose another assumption to bound the difference between  $\widetilde{H}_{\mu,i}(x_i)$ and $\nabla^2 f_i(x_i)$.
\begin{assumption}\label{assump2}
	The estimated Hessian $\widetilde{H}_{\mu,i}(x_i)$ satisfies
	\vspace{-0.15cm}\begin{equation*}
		\theta\widetilde{H}_{\mu,i}(x_i)\preceq\nabla^2f_i(x_i)\preceq(2-\theta)\widetilde{H}_{\mu,i}(x_i),
		\vspace{-0.15cm}\end{equation*}
	for $i=1,\dots,N$ and some $\theta\in(0,1]$.
\end{assumption}

Parameter $\theta$ measures how accurate $\widetilde{H}_{\mu,i}(x_i)$ approximates $\nabla^2f_i(x_i)$. Specifically, $\widetilde{H}_{\mu,i}(x_i)$ reduces to the exact Hessian $\nabla^2f_i(x_i)$ when $\theta=1$. 
The way of constructing zeroth-order estimate for Hessian (\ref{hessian_estimate_i}) may satisfy Assumption \ref{assump2} with proper parameter values such as sufficiently large $b$, small $\mu$ and evenly distributed $u_j$.

From Assumption \ref{assump1} and \ref{assump2}, we have 
\vspace{-0.12cm}\begin{equation}
	\frac{1}{2-\theta}\Lambda_m\preceq\frac{1}{2-\theta}\nabla^2 f(\mathbf{x})\preceq\widetilde{H}_\mu(\mathbf{x})\preceq\frac{1}{\theta}\nabla^2 f(\mathbf{x})\preceq\frac{1}{\theta}\Lambda_M. \label{H_range}
\vspace{-0.12cm}\end{equation}
Let $\bar{\Lambda}=\underline{\alpha}^{-1}\left(\frac{1}{\theta}\Lambda_M-\frac{\Lambda_M+\Lambda_m}{2}\right)$. The convergence analysis relies on the following condition
\vspace{-0.12cm}\begin{align}
	D\succ&\frac{\Lambda_M}{2\eta}+\rho(W+\mathbf{I}_{Nd})+\left(\frac{2}{\theta}-\frac{3}{2}\right)\Lambda_M-\frac{3}{2}\Lambda_m \nonumber\\
	&+\left(\frac{1}{\theta}\Lambda_M-\frac{\Lambda_M+\Lambda_m}{2}\right)^2 \label{eq_lemma3_1},
\vspace{-0.12cm}\end{align}
for any $\eta>1$. With (\ref{eq_lemma3_1}), it is guaranteed that $\widetilde{H}_{\mu}(\mathbf{x})+D\succ\mathbf{O}_{Nd}\ \forall \mathbf{x}\in\mathbb{R}^{Nd}$ since $\widetilde{H}_{\mu}(\mathbf{x})\succeq\frac{1}{2-\theta}\Lambda_m$.

For better presentation, let $\mathbf{z}^k=((\mathbf{x}^k)^T,(\mathbf{v}^k)^T)^T$ and $\mathbf{z}^*=((\mathbf{x}^*)^T,(\mathbf{v}^*)^T)^T$. Also, let $\lambda_W>0$ be the smallest nonzero eigenvalue of $W$. The main convergence result of ZoPro is provided below.

\begin{figure*}[!htbp]
	\centering
	\includegraphics[width=0.84\textwidth,height=0.21\textheight]{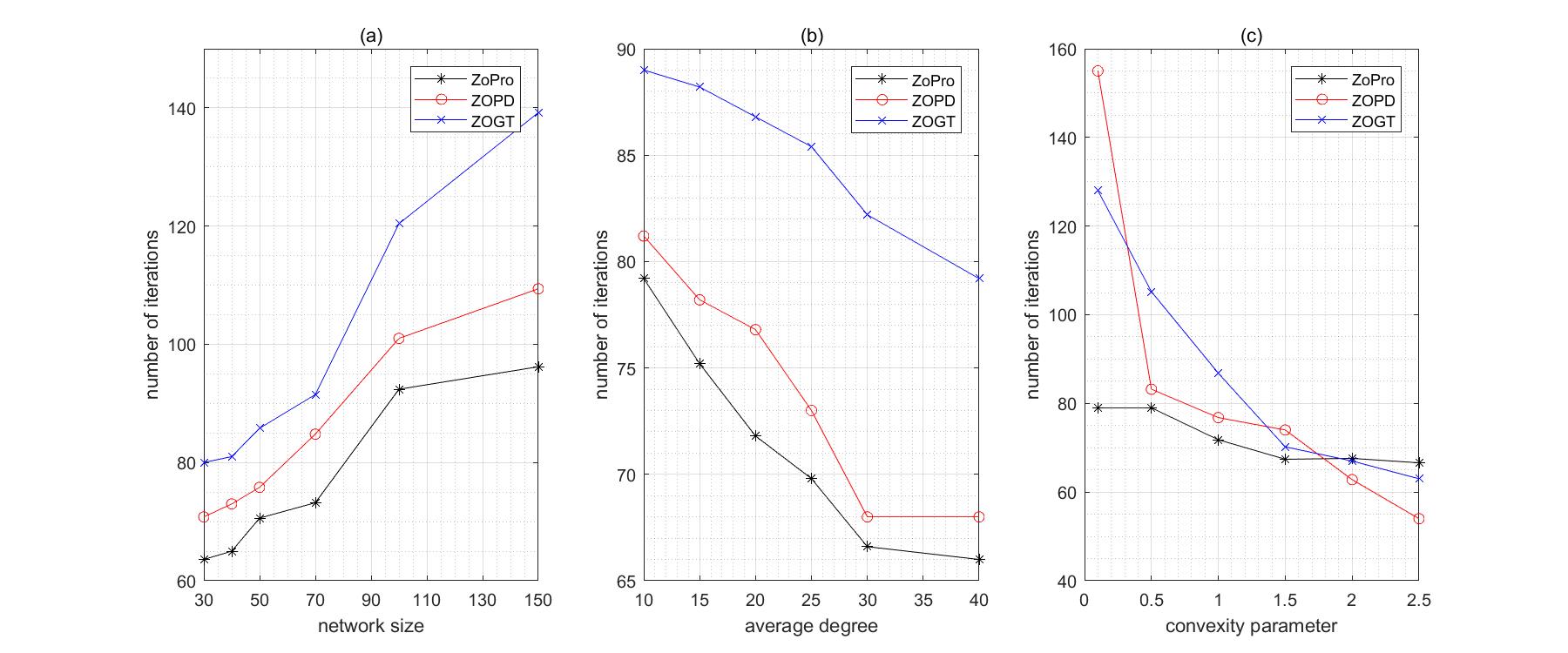}
	\caption{Convergence performance of ZoPro, ZOPD and ZOGT}
	\label{Convergence performance of ZoPro, ZOPD and ZOGT}
\end{figure*}

\begin{theorem}\label{theorem1}
	 Suppose Assumptions \ref{assump1} and \ref{assump2} hold. Assume (\ref{eq_lemma3_1}) holds for some $\eta>1$. Then, for any $\beta>\underline{\alpha}^{-1}$ and $\gamma>\frac{2m(\eta-1)+\eta+\beta}{\eta-1}$, $\mathbf{z}^k$ converges linearly to a neighborhood of $\mathbf{z}^*$ in expectation, i.e. there exists $\delta\in(0,1)$ such that for each $k\geq 0$,
	\vspace{-0.15cm}\begin{align}
		&\mathbf{E}\left[\left\|\mathbf{z}^{k+1}-\mathbf{z}^*\right\|_Q^2\right]\leq(1-\delta)\mathbf{E}\left[\left\|\mathbf{z}^k-\mathbf{z}^*\right\|_Q^2\right]+G, \label{theorem1_1}\\
		&\lim\sup_{k\rightarrow\infty}\mathbf{E}\left[\left\|\mathbf{z}^k-\mathbf{z}^*\right\|_Q^2\right]\leq\frac{G}{\delta}. \label{theorem1_7}
	\end{align}
	In particular, given any $c_1>0$,  $G=\rho(\eta+\frac{1-\eta}{\gamma})G_2^2+2(G_1^2+G_2^2)+\frac{2\delta(1+c_1)(G_1+G_2)^2}{\lambda_W}$ and
	\begin{align}
		\delta&=\mathop{\sup}\limits_{c_1,c_2>0}\min\left\{\frac{\rho\lambda_W\kappa_{\beta,\eta}}{2\underline{\alpha}^{-2}(1+c_1)\left\|\frac{1}{\theta}\Lambda_M+D\right\|^2} \right., \notag \\
		&\left.{\frac{1}{(1+1/c_1)(1+c_2)},\frac{\delta_c}{{\lambda_{\max}(\mathcal{B}/\rho)}}}\right\} \label{theorem1_2},	    
	\end{align}
	in which $\mathcal{B}=\frac{(1+1/c_1)(1+1/c_2)\Lambda_M^2}{\lambda_W}+\rho R$, $\delta_c=(2m-\gamma)(1-\eta)-\eta-\beta$ and $\kappa_{\beta,\eta}=\lambda_{\min}(R-\frac{\Lambda _M}{2\eta}-\frac{\bar{\Lambda}^2}{\beta}- 2\bar{\Lambda}-\rho(\mathbf{I}_{Nd}+W))>0$.
\end{theorem} 

\begin{proof}
	See Appendix \ref{Proof of Theorem1}.
\end{proof}

Subsequently, we discuss the influence of the objective function and the network topology on the convergence rate of ZoPro. From (\ref{theorem1_1}), note that $\delta$ mainly depends on $M$, $m$ and $\lambda_{W}$. To see this, let $f_i$ $\forall\ i\in\mathcal{V}$ be identically strongly convex with parameter $m$ and smooth with parameter $M$ such that $0<m<M$. Let $W=\mathbf{I}_{Nd}-A\otimes \mathbf{I}_d$, where $A=A^T$ is a doubly stochastic matrix. It can be shown that larger $m$, smaller $M$ and larger $\lambda_{W}$ (which suggests denser connectivity of $\mathcal{G}$) lead to larger $\delta$ and a faster convergence speed of ZoPro.

Also, we discuss the factors that affect the ultimate optimality error, i.e., the expected distance between $\mathbf{z}^k$ and $\mathbf{z}^*$ as $k\rightarrow\infty$. From (\ref{theorem1_7}), this expected distance mainly depends on $\mu$, $b$ and $M$. It can be shown that smaller $\mu$ (accurate smoothed approximation), larger $b$ (enough sample points for zeroth-order oracle) and smaller $M$ (well-conditioned objective function) contribute to a smaller expected error. However, we need to control $b$ within a moderate range to avoid high computational cost in constructing Hessian and gradient estimates.

\section{NUMERICAL EXPERIMENTS}
This section illustrates the practical convergence performance of ZoPro and its comparisons with related algorithms.

In the numerical experiment, we consider the following logistic regression problem with $\mathcal{L}_2$ regularization: All the nodes need to cooperatively minimize the objective function $f(x)=\sum_{i\in\mathcal{V}}\left(\frac{\lambda}{2N}\|x_i\|^2+\sum_{l=1}^{q_i}\text{log}\left(1+\text{exp}(-v_{il}\mathbf{u}_{il}^Tx_i)\right)\right)$, where $x_i\in\mathbb{R}^d$, $\lambda$ is a regularized parameter, $N$ is the number of nodes, $q_i$ is the sample number assigned to node $i$, $\mathbf{u}_{ij}$, $j=1,\dots,q_i$ is the data and $v_{ij}\in\{-1,+1\}$, $j=1,\dots,q_i$ is the corresponding label. Here, we set $q_i=5\ \forall i\in\mathcal{V}$.
\vspace{-0.1cm}

\subsection{Comparison with Zeroth-order Methods}
We compare the convergence performance of ZoPro with some other existing zeroth-order optimization methods, including zeroth-order gradient-tracking method (ZOGT) \cite{b22} and zeroth-order primal-dual method (ZOPD) \cite{b23}.

We set $d=20$. In order to test how the problem and network characteristics influence the convergence performance, we consider three parameters, i.e., the network size $N$, the average node degree $d_a=\sum_{i\in\mathcal{V}}|\mathcal{N}_i|/N$ and the convexity parameter $\lambda$. Accordingly, we run three groups of experiments, each of which fixes two of these parameters and varies the other. We express each experiment as a triplet $(N,d_a,\lambda)$, and set the three experiment groups to (G1) $(N,20,1)$, $N=30,40,50,70,100,150$, (G2) $(50,d_a,1)$, $d_a=10,15,20,25,30,40$ and (G3) $(50,20,\lambda)$, $\lambda=0.1,0.5,1,1.5,2,2.5$. For each value of $(N,d_a,\lambda)$, we generate 10 random scenarios and plot their average in our figures. In each scenario, the undirected network is randomly generated with the given $N$ and $d_a$. 

The algorithm parameters of ZoPro are selected moderately. We let the smoothness parameter $\mu=0.05$, the batch size $b=50$ and the stepsize control parameter $c=0.1$. We terminate the algorithms when the average error of all the nodes $\sum_{i\in\mathcal{V}}\|x_i^k-x^*\|^2/N$ drops below $10^{-4}$ and remains there for $100$ more iterations. Therefore, we define the number of iterations needed for convergence as $\min\left\{k:\sum_{i\in\mathcal{V}}\|x_i^{k+t}-x^*\|^2/N\leq 10^{-4},\ 0\le t\le 100\right\}$.

Figure \ref{Convergence performance of ZoPro, ZOPD and ZOGT} plots the number of iterations needed for convergence of ZoPro, ZOPD and ZOGT with $(N,d_a,\lambda)$ given by (G1), (G2) and (G3), respectively. Observe that smaller and denser networks as well as larger convexity parameters essentially lead to faster convergence for ZoPro, ZOGT and ZOPD. Compared to ZOGT and ZOPD, our proposed ZoPro requires the fewest iterations to reach the convergence criterion in most of the cases.

\subsection{Comparison with Second-order Methods}

To illustrate the computational efficiency of ZoPro, we make a comparison of ZoPro and the well-noted second-order methods ESOM \cite{b17}, DQM\cite{b18} and SoPro\cite{b20} in terms of running time. Here, we consider the same problem form as before, and set $(N,d_a,\lambda,d)$ to be $(30,10,1,20)$ and $(100,20,1,30)$ to simulate a medium-scale problem and a large-scale problem. 

\begin{figure}[!htpb]
	\centering
	\includegraphics[width=0.36\textwidth,height=0.163\textheight]{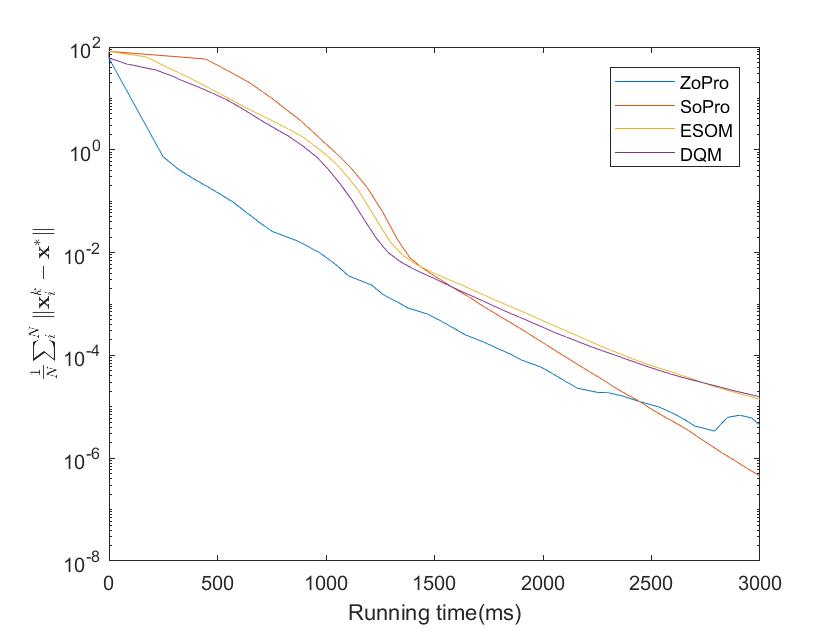}
	\caption{Convergence performance of ZoPro, ESOM, DQM and SoPro for the medium-scale problem}
	\label{Convergence performance of ZoPro, ESOM, DQM and SoPro in medium-scaled problem}
	\vspace{-0.15cm}
\end{figure}
\vspace{-0.15cm}

\begin{figure}[!htpb]
	\centering
	\includegraphics[width=0.36\textwidth,height=0.168\textheight]{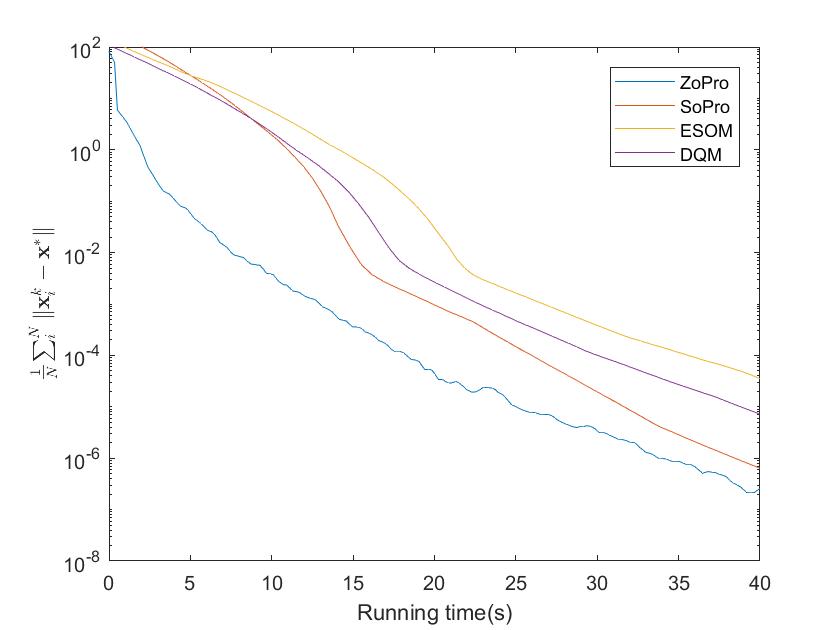}
	\caption{Convergence performance of ZoPro, ESOM, DQM and SoPro for the large-scale problem}
	\label{Convergence performance of ZoPro, ESOM, DQM and SoPro in larger-scaled problem}
\end{figure}

In Figure~\ref{Convergence performance of ZoPro, ESOM, DQM and SoPro in medium-scaled problem} for the medium-scale problem, ZoPro consistently takes shorter running time for reaching the same accuracy than ESOM and DQM. Compared to SoPro which ZoPro originates from, ZoPro requires less time to reach the accuracy $10^{-4}$. However, when the accuracy is enhanced to $10^{-6}$, SoPro is a better option. In Figure~\ref{Convergence performance of ZoPro, ESOM, DQM and SoPro in larger-scaled problem} for the large-scale problem, the advantage of ZoPro against all the other three second-order methods becomes more prominent and ZoPro takes the shortest time in achieving any accuracy above $10^{-6}$.

\section{CONCLUSION}

We develop a zeroth-order proximal algorithm (ZoPro) for solving consensus optimization problems over undirected networks. ZoPro approximates exact gradients and Hessian matrices in a powerful second-order method SoPro using a zeroth-order oracle, which significantly reduces the computational complexity, particularly in solving large-scale problems. ZoPro inherits some appealing features of SoPro, including full decentralization and fast convergence. We show that ZoPro achieves linear convergence to a neighborhood of the optimum in expectation when the problem is strongly convex and smooth. We also demonstrate the fast convergence and computational efficiency of ZoPro by comparing it with several state-of-the-art distributed zeroth-order and second-order algorithms through extensive simulations.

\section{APPENDIX} \label{appendix}
\subsection{Proof of Proposition~\ref{prop1}}\label{Proof of Proposition1}
We assume that none of (i), (ii) and (iii) hold and establish a contradiction. Since (i) does not occur, the decreasing sequence $\left\{f_i(x_i^k)\right\}$ is bounded below and we have $f_i(x_i^k)\rightarrow \bar{f_i}$ as $k\rightarrow \infty$ for some $\bar{f_i}$. In particular, $f_i(x_i^{k+1})-f_i(x_i^k)\rightarrow 0$ as $k\rightarrow \infty$. Next, since (ii) and (iii) do not occur, there are a subsequence $J\subseteq \mathbb{N}$ and a vector $\bar{d_i}$ such that $d_i^{k}\rightarrow \bar{d_i}$ for $k\in J$ and 
\vspace{-0.2cm}\begin{equation*}
		\sup_{k\in J}f_i'(x_i^k;d_i^k)=:\beta_0<0.
		\vspace{-0.2cm}\end{equation*}
The Armijo condition and the fact that $f_i(x_i^{k+1})-f_i(x_i^k)\rightarrow 0$ as
 $k\rightarrow \infty$ imply that
\vspace{-0.2cm}\begin{equation*}
		\alpha_i^{k}f_i'(x_i^k;d_i^k)\rightarrow 0.
		\vspace{-0.2cm}\end{equation*}
	Since $f_i'(x_i^k;d_i^k)\leq\beta_0<0$ for $k\in J$, we have $\alpha_i^k\rightarrow 0$ as $k\rightarrow \infty$. With no loss in generality, we assume that $\alpha_i^k\in(0,1)$ for all $k\in J$. For $\epsilon\in(0,1)$, we have
	\vspace{-0.2cm}\begin{equation}
		\epsilon\gamma_0^{-1}\alpha_i^k f_i'(x_i^k;d_i^k)<f_i(x_i^{k}+\alpha_i^{k}\gamma_0^{-1}d_i^{k})-f_i(x_i^k) \label{armijo}
		\vspace{-0.15cm}\end{equation}
	for all $k\in J$, where $\gamma_0>0$ is the search control parameter. By the Mean Value Theorem, there exists $\theta_0\in (0,1)$ for each $k\in J$ that
	\vspace{-0.15cm}\begin{equation*}
		f_i(x_i^k+\alpha_i^k\gamma_0^{-1}d_i^k)-f_i(x_i^k)=\alpha_i^k\gamma_0^{-1}f_i'(x_i^k+\theta_0\alpha_i^k\gamma_0^{-1}d_i^k;d_i^k).
		\vspace{-0.15cm}\end{equation*}
	For better presentation, we define $\hat{x}_i^k=x_i^k+\theta_0\alpha_i^k\gamma_0^{-1}d_i^k$. Now, since $f_i$ is L-smooth, we have
	\vspace{-0.15cm}\begin{align}
		&\alpha_i^k\gamma_0^{-1}f_i'(\hat{x}_i^k;d_i^k) \nonumber\\
		=&\alpha_i^k\gamma_0^{-1}f_i'(x_i^k;d_i^k)+\alpha_i^k\gamma_0^{-1}\left(f_i'(\hat{x}_i^k;d_i^k)-f_i'(x_i^k;d_i^k)\right) \nonumber\\
		=&\alpha_i^k\gamma_0^{-1}f_i'(x_i^k;d_i^k)+\alpha_i^k\gamma_0^{-1}\left(\nabla f_i(\hat{x}_i^k)-\nabla f_i(x_i^k)\right)^Td_i^k \nonumber\\
		\leq&\alpha_i^k\gamma_0^{-1}f_i'(x_i^k;d_i^k)+\alpha_i^k\gamma_0^{-1}L\|\hat{x}_i^k-x_i^k\|\|d_i^k\| \nonumber\\
		=&\alpha_i^k\gamma_0^{-1}f_i'(x_i^k;d_i^k)+(\alpha_i^k\gamma_0^{-1})^2L\theta_0\|d_i^k\|^2. \label{armijo_2}
		\vspace{-0.15cm}\end{align}
	Combining (\ref{armijo}) and (\ref{armijo_2}), we have
	\vspace{-0.15cm}\begin{equation*}
		\epsilon\gamma_0^{-1}\alpha_i^kf_i'(x_i^k;d_i^k)<\alpha_i^k\gamma_0^{-1}f_i'(x_i^k;d_i^k)+(\alpha_i^k\gamma_0^{-1})^2L\theta_0\|d_i^k\|^2,
		\vspace{-0.15cm}\end{equation*}
	which implies that
	\vspace{-0.15cm}\begin{equation*}
		0<(1-\epsilon)\beta_0+(\alpha_i^k\gamma_0^{-1})L\theta_0\|d_i^k\|^2,\ \forall k\in J.
		\vspace{-0.15cm}\end{equation*}
	Taking the limit over $k\in J$, we obtain the contradiction
	\vspace{-0.15cm}\begin{equation*}
		0\leq(1-\epsilon)\beta_0<0,
		\vspace{-0.15cm}\end{equation*}
	where the first inequality is due to the fact that $\alpha_i^k\rightarrow 0$ as $k\rightarrow\infty$ and the second inequality results from $\epsilon\in(0,1)$ and $\beta_0<0$.

\subsection{Proof of Theorem~\ref{theorem1}}\label{Proof of Theorem1}
Similar to \cite{b20}, since the dual optimum of problem (\ref{consensus_problem}) can be arbitrarily chosen such that $\nabla f(\mathbf{x}^*)=-W^{\frac{1}{2}}\mathbf{v}$, we simply define
\vspace{-0.15cm}\begin{equation}\label{dual_optimum}
	\mathbf{v}^*=-(W^\dagger)^{\frac{1}{2}}\nabla f(\mathbf{x}^*)
	\vspace{-0.15cm}\end{equation} 
as a particular dual optimum. Besides, let $\mathbf{v}^k=(W^{\dagger})^{\frac{1}{2}}\mathbf{q}^k$ and we have 
\vspace{-0.15cm}\begin{equation}\label{v_on}
	\mathbf{v}^*,\mathbf{v}^k,\mathbf{v}^k-\mathbf{v}^*\in\left\{\mathbf{x}\in\mathbb{R}^{Nd}|x_1+\dots+x_n=\mathbf{0}_d\right\}.
	\vspace{-0.15cm}\end{equation}

First, we derive a lemma to bound the difference between $\mathbf{E}\left[\|\mathbf{z}^k-\mathbf{z}^*\|^2_Q\right]$ and $\mathbf{E}\left[\|\mathbf{z}^{k+1}-\mathbf{z}^*\|^2_Q\right]$.

\begin{lemma}\label{lemma3}
	Suppose Assumptions \ref{assump1} and \ref{assump2} hold. For each $k\geq0$ and any $\eta,\beta,\gamma>0$, we have
	\vspace{-0.15cm}\begin{align}
		&\mathbf{E}\left[\left\|\mathbf{z}^k-\mathbf{z}^*\right\|_Q^2\right]-\mathbf{E}\left[\left\|\mathbf{z}^{k+1}-\mathbf{z}^*\right\|_Q^2\right] \notag\\
		\geq&\rho\delta_c\mathbf{E}\left[\left\|\mathbf{x}^k-\mathbf{x}^*\right\|^2\right]+\rho^2\mathbf{E}\left[\left\|\mathbf{x}^k\right\|_W^2\right]-\rho(\eta+\frac{1-\eta}{\gamma})G_2^2 \notag\\
		&-\rho\mathbf{E}\left[\left\|\mathbf{x}^{k+1}-\mathbf{x}^k\right\|_{\mathcal{A}_{\beta,\eta}+\rho W-R}^2\right]-2(G_1^2+G_2^2),\label{eq_lemma3_2}
		\vspace{-0.15cm}\end{align}	
	where $\mathcal{A}_{\beta,\eta}=\frac{\Lambda_M}{2\eta}+\frac{\bar{\Lambda}^2}{\beta}+2\bar{\Lambda}+\rho\mathbf{I}_{Nd}$.
\end{lemma}

\begin{proof}
	From (\ref{zeroth_dual_iteration}) and   $\mathbf{q}^k=W^{\frac{1}{2}}\mathbf{v}^k$, we have
	\vspace{-0.15cm}\begin{equation}
		\mathbf{v}^{k+1}=\mathbf{v}^{k}+\rho W^{\frac{1}{2}}\mathbf{x}^{k+1}.
		\label{eq_dual_iteration}
		\vspace{-0.15cm}\end{equation}
	Using (\ref{eq_dual_iteration}) and $W\mathbf{x}^*=\textbf{0}_{Nd}$, we can derive
	\vspace{-0.15cm}\begin{align}
		&\left\langle\mathbf{v}^k-\mathbf{v}^{k+1},\mathbf{v}^{k+1}-\mathbf{v}^*\right\rangle \notag\\
		=&-\rho\left\langle\mathbf{x}^{k+1}-\mathbf{x}^*,W^{\frac{1}{2}}(\mathbf{v}^{k+1}-\mathbf{v}^*)\right\rangle.\label{eq_lemma2_1}
		\vspace{-0.15cm}\end{align}
	From (\ref{zeroth_dual_iteration}) and (\ref{zeroth_primal_backtracking_iteration}),
	\vspace{-0.15cm}\begin{align}
		&W^{\frac{1}{2}}\mathbf{v}^{k+1}=W^{\frac{1}{2}}(\mathbf{v}^{k+1}-\mathbf{v}^k)+W^{\frac{1}{2}}\mathbf{v}^k \notag\\
		=&\left(\rho W-(A^k)^{-1}H^k\right)(\mathbf{x}^{k+1}-{\mathbf{x}^k})-\widetilde{g}_\mu({\mathbf{x}^k}). \label{eq_lemma2_2}
		\vspace{-0.15cm}\end{align}
	Combining (\ref{dual_optimum}), (\ref{eq_lemma2_1}) and (\ref{eq_lemma2_2}), we have
	\vspace{-0.15cm}\begin{align}
		&\left\langle\mathbf{v}^k-\mathbf{v}^{k+1},\mathbf{v}^{k+1}-\mathbf{v}^*\right\rangle \notag\\
		=&\rho\left\langle\mathbf{x}^{k+1}-\mathbf{x}^*,\widetilde{g}_\mu(\mathbf{x}^k)-\nabla f(\mathbf{x}^*)\right\rangle \notag\\
		&+\rho\left\langle\mathbf{x}^{k+1}-\mathbf{x}^*,\left((A^k)^{-1}H^k-\rho W\right)(\mathbf{x}^{k+1}-\mathbf{x}^k)\right\rangle. \label{eq_lemma2_3}
		\vspace{-0.15cm}\end{align}
	Moreover, since $W\mathbf{x}^*=0_{Nd}$, it follows that
	\vspace{-0.15cm}\begin{align}
		&- \left\langle \mathbf{x}^{k+1}-\mathbf{x}^*,W(\mathbf{x}^{k+1}-\mathbf{x}^k)\right\rangle+\left\|\mathbf{x}^{k+1}\right\|_W^2 \notag\\
		=&\frac{1}{2}(\left\|\mathbf{x}^{k+1}\right\|_W^2+\left\|\mathbf{x}^k\right\|_W^2 -\left\|\mathbf{x}^{k+1}-\mathbf{x}^k \right\|_W^2).\label{eq_lemma2_4}
		\vspace{-0.15cm}\end{align}
	We expand the left-hand side of (\ref{eq_lemma3_2}) and obtain
	\vspace{-0.15cm}\begin{align}
		&\mathbf{E}\left[\left\|\mathbf{z}^k-\mathbf{z}^*\right\|_Q^2\right]-\mathbf{E}\left[\left\|\mathbf{z}^{k+1}-\mathbf{z}^*\right\|_Q^2\right] \nonumber\\
		=&\mathbf{E}\left[\left\|\mathbf{z}^k-\mathbf{z}^{k+1}\right\|_Q^2\right]+2\mathbf{E}\left[\rho\left\langle \mathbf{x}^{k+1}-\mathbf{x}^*,R(\mathbf{x}^k-\mathbf{x}^{k+1})\right\rangle\right] \nonumber\\
		+&2\mathbf{E}\left[\left\langle\mathbf{v}^k-\mathbf{v}^{k+1},\mathbf{v}^{k+1}-\mathbf{v}^*\right\rangle\right]. \label{eq_lemma2_5}
	\end{align}
	By incorporating (\ref{eq_lemma2_4}) into (\ref{eq_lemma2_3}) and combining the resulting equation with (\ref{eq_lemma2_5}), we have
	\begin{align}
		&\mathbf{E}\left[\left\|\mathbf{z}^k-\mathbf{z}^*\right\|_Q^2\right]-\mathbf{E}\left[\left\|\mathbf{z}^{k+1}-\mathbf{z}^*\right\|_Q^2\right] \notag\\
		=&2\rho\mathbf{E}\left[\left\langle\mathbf{x}^{k+1}-\mathbf{x}^*,\widetilde{g}_\mu(\mathbf{x}^k)-\nabla f(\mathbf{x}^*)\right\rangle\right]+\mathbf{E}\left[\rho^2\left\|\mathbf{x}^k\right\|_W^2\right] \notag\\
		&+2\rho\mathbf{E}\left[\left\langle\mathbf{x}^{k+1}-\mathbf{x}^*,\left((A^k)^{-1}H^k-R\right)(\mathbf{x}^{k+1}-{\mathbf{x}^k})\right\rangle\right] \notag\\
		&+\rho\mathbf{E}\left[\left\|\mathbf{x}^{k+1}-\mathbf{x}^k\right\|_{R-\rho W}^2\right].  \label{eq_lemma2_6}
		\vspace{-0.15cm}\end{align}
	Now, we need to bound the first and third terms of the right-hand side of (\ref{eq_lemma2_6}). First, based on the AM-GM inequality, (\ref{norm_bound}) and (\ref{norm_bound_2}), we provide a lower bound on the first term. For any $\eta>0$, we have
	\vspace{-0.15cm}\begin{align}
		&\mathbf{E}\left[\left\langle\mathbf{x}^{k+1}-\mathbf{x}^k,\widetilde{g}_\mu(\mathbf{x}^k)-\nabla f(\mathbf{x}^*)\right\rangle\right] \nonumber\\
		=&\mathbf{E}\left[\left\langle\mathbf{x}^{k+1}-\mathbf{x}^k,\nabla f(\mathbf{x}^k)-\nabla f(\mathbf{x}^*)\right\rangle\right] \nonumber\\
		&+\mathbf{E}\left[\left\langle\mathbf{x}^{k+1}-\mathbf{x}^k,\widetilde{g}_\mu(\mathbf{x}^k) 
		-\nabla f_\mu(\mathbf{x}^k)\right\rangle\right] \nonumber\\ 
		&+\mathbf{E}\left[\left\langle\mathbf{x}^{k+1}-\mathbf{x}^k,\nabla f_\mu(\mathbf{x}^k)-\nabla f(\mathbf{x}^k)\right\rangle\right] \nonumber\\
		\geq&-\eta\mathbf{E}\left[\|\nabla f(\mathbf{x}^k)-\nabla f(\mathbf{x}^*)\|_{\Lambda_{M^{-1}}}^2\right] \nonumber\\
		&-\frac{\mathbf{E}\left[\|\mathbf{x}^{k+1}-\mathbf{x}^k\|_{\Lambda_M}^2\right]}{4\eta} \nonumber\\
		&-\frac{1}{\rho}\mathbf{E}\left[\left\|\widetilde{g}_\mu(\mathbf{x}^k)-\nabla f_\mu(\mathbf{x}^k)\right\|^2\right]-\frac{\rho}{4}\mathbf{E}\left[\left\|\mathbf{x}^{k+1}-\mathbf{x}^k\right\|^2\right] \nonumber\\
		&-\frac{1}{\rho}\mathbf{E}\left[\left\|\nabla f_\mu(\mathbf{x}^k)-\nabla f(\mathbf{x}^k)\right\|^2\right]-\frac{\rho}{4}\mathbf{E}\left[\left\|\mathbf{x}^{k+1}-\mathbf{x}^k\right\|^2\right] \nonumber\\
		\geq&-\eta\mathbf{E}\left[\|\nabla f(\mathbf{x}^k)-\nabla f(\mathbf{x}^*)\|_{\Lambda_{M^{-1}}}^2\right] \nonumber\\
		&-\mathbf{E}\left[\left\|\mathbf{x}^{k+1}-\mathbf{x}^k\right\|^2_{\frac{\Lambda_M}{4\eta}+\frac{\rho \mathbf{I}_{Nd}}{2}}\right]-\frac{G_1^2+G_2^2}{\rho}. \label{eq_lemma2_9}
		\vspace{-0.15cm}\end{align}
	Besides, due to the gradient Lipschitz of $f$ and unbiasedness of $\widetilde{g}_\mu$, we have
	\vspace{-0.1cm}\begin{align*}
		&\mathbf{E}\left[\left\langle\mathbf{x}^k-\mathbf{x}^*,\widetilde{g}_\mu(\mathbf{x}^k)-\nabla f(\mathbf{x}^*)\right\rangle\right] \nonumber \\
		&+\mathbf{E}\left[\left\langle\mathbf{x}^k-\mathbf{x}^*,\nabla f(\mathbf{x}^k)-\nabla f(\mathbf{x}^*)\right\rangle\right] \nonumber\\
		\geq&\mathbf{E}\left[\|\nabla f(\mathbf{x}^k)-\nabla f(\mathbf{x}^*)\|_{\Lambda_{M^{-1}}}^2\right] \nonumber\\
		&-\frac{1}{2}\mathbf{E}\left[\left\|\mathbf{x}^k-\mathbf{x}^*\right\|^2\right]-\frac{1}{2}\mathbf{E}\left[\left\|\nabla f_\mu(\mathbf{x}^k)-\nabla f(\mathbf{x}^k)\right\|^2\right] \nonumber\\
		\geq&\mathbf{E}\left[\|\nabla f(\mathbf{x}^k)-\nabla f(\mathbf{x}^*)\|_{\Lambda_{M^{-1}}}^2\right] \nonumber\\
		&-\frac{1}{2}\mathbf{E}\left[\left\|\mathbf{x}^k-\mathbf{x}^*\right\|^2\right]-\frac{1}{2}G_2^2. 
	\end{align*}
	Multiply this inequality by $\eta$ and add it to (\ref{eq_lemma2_9}). For any $\eta>0$ and $\gamma>0$, we then have
	\begin{align}
		&\mathbf{E}\left[\left\langle\mathbf{x}^{k+1}-\mathbf{x}^*,\widetilde{g}_{\mu}(\mathbf{x}^k)-\nabla f(\mathbf{x}^*)\right\rangle\right] \notag\\
		\geq&(1-\eta)\mathbf{E}\left[\left\langle\mathbf{x}^k-\mathbf{x}^*,\nabla f(\mathbf{x}^k)-\nabla f(\mathbf{x}^*)\right\rangle\right] \nonumber\\
		&+(1-\eta)\mathbf{E}\left[\left\langle\mathbf{x}^k-\mathbf{x}^*,\nabla f_\mu(\mathbf{x}^k)-\nabla f(\mathbf{x}^k)\right\rangle\right] \nonumber\\
		&-\mathbf{E}\left[\left\|\mathbf{x}^{k+1}-\mathbf{x}^k\right\|^2_{\frac{\Lambda_M}{4\eta}+\frac{\rho \mathbf{I}_{Nd}}{2}}\right] \nonumber\\
		&-\frac{1}{2}\eta\mathbf{E}\left[\left\|\mathbf{x}^k-\mathbf{x}^*\right\|^2\right]-\frac{1}{2}\eta G_2^2-\frac{G_1^2+G_2^2}{\rho} \nonumber\\
		\geq&(1-\eta)\mathbf{E}\left[\left\langle\mathbf{x}^k-\mathbf{x}^*,\nabla f(\mathbf{x}^k)-\nabla f(\mathbf{x}^*)\right\rangle\right] \nonumber\\
		&-\frac{1}{2}(1-\eta)\gamma\mathbf{E}\left[\left\|\mathbf{x}^k-\mathbf{x}^*\right\|^2\right] \nonumber\\
		&-\frac{1}{2}(1-\eta)\frac{1}{\gamma}\mathbf{E}\left[\left\|\nabla f_\mu(\mathbf{x}^k)-\nabla f(\mathbf{x}^k)\right\|^2\right] \nonumber\\
		&-\mathbf{E}\left[\left\|\mathbf{x}^{k+1}-\mathbf{x}^k\right\|^2_{\frac{\Lambda_M}{4\eta}+\frac{\rho \mathbf{I}_{Nd}}{2}}\right] \nonumber\\
		&-\frac{1}{2}\eta\mathbf{E}\left[\left\|\mathbf{x}^k-\mathbf{x}^*\right\|^2\right]-\frac{1}{2}\eta G_2^2-\frac{G_1^2+G_2^2}{\rho} \nonumber\\
		\geq&(1-\eta)\mathbf{E}\left[\left\langle\mathbf{x}^k-\mathbf{x}^*,\nabla f(\mathbf{x}^k)-\nabla f(\mathbf{x}^*)\right\rangle\right] \nonumber\\
		&-\mathbf{E}\left[\left\|\mathbf{x}^{k+1}-\mathbf{x}^k\right\|^2_{\frac{\Lambda_M}{4\eta}+\frac{\rho \mathbf{I}_{Nd}}{2}}\right] \nonumber\\
		&-\frac{1}{2}\left((1-\eta)\gamma+\eta\right)\mathbf{E}\left[\left\|\mathbf{x}^k-\mathbf{x}^*\right\|^2\right] \nonumber\\
		&-\frac{1}{2}\left(\eta+\frac{1-\eta}{\gamma}\right)G_2^2-\frac{G_1^2+G_2^2}{\rho}.
		\label{eq_lemma2_7}
	\end{align}
	In addition, we have $(A^k)^{-1}H^k-R=(A^k)^{-1}\widetilde{H}_\mu(\mathbf{x}^k)+D\left((A^k)^{-1}-\underline{\alpha}^{-1}\mathbf{I}_{Nd}\right)-\underline{\alpha}^{-1}\frac{\Lambda_m+\Lambda_M}{2}\preceq\underline{\alpha}^{-1}\widetilde{H}_\mu(\mathbf{x}^k)-\underline{\alpha}^{-1}\frac{\Lambda_m+\Lambda_M}{2}$. Due to (\ref{H_range}), $(A^k)^{-1}H^k-R\preceq\bar{\Lambda}$. For any $\beta>0$, we obtain
	\vspace{-0.1cm}
	\begin{align}
		&\mathbf{E}\left[\left\langle\mathbf{x}^{k+1}-\mathbf{x}^*,\left((A^k)^{-1}H^k-R\right)(\mathbf{x}^{k+1}-\mathbf{x}^k)\right\rangle\right] \notag \\
		=&\mathbf{E}\left[\left\langle\mathbf{x}^k-\mathbf{x}^*,\left((A^k)^{-1}H^k-R\right)(\mathbf{x}^{k+1}-\mathbf{x}^k)\right\rangle\right] \notag \\
		&+\mathbf{E}\left[\left\langle\mathbf{x}^{k+1}-\mathbf{x}^k,\left((A^k)^{-1}H^k-R\right)(\mathbf{x}^{k+1}-\mathbf{x}^k)\right\rangle\right] \notag \\
		\geq&-\frac{\beta}{2}\mathbf{E}\left[\|\mathbf{x}^k-\mathbf{x}^*\|^2\right]-\frac{1}{2\beta}\mathbf{E}\left[\|\mathbf{x}^{k+1}-\mathbf{x}^k\|_{\bar{\Lambda}^2}^2\right] \notag\\
		&-\mathbf{E}\left[\|\mathbf{x}^{k+1}-\mathbf{x}^k\|_{\bar{\Lambda}}^2\right]. \label{eq_lemma2_8}
		\vspace{-0.1cm}
	\end{align}
	Combining (\ref{eq_lemma2_6}), (\ref{eq_lemma2_7}) and (\ref{eq_lemma2_8}) yields the following inequality
	\vspace{-0.2cm}\begin{align}
		&\mathbf{E}\left[\left\|\mathbf{z}^k-\mathbf{z}^*\right\|_Q^2\right]-\mathbf{E}\left[\left\|\mathbf{z}^{k+1}-\mathbf{z}^*\right\|_Q^2\right] \notag\\
		\geq&2\rho(1-\eta)\mathbf{E}\left[\left\langle\mathbf{x}^k-\mathbf{x}^*,\nabla f(\mathbf{x}^k)-\nabla f(\mathbf{x}^*)\right\rangle\right] \notag\\
		&+\rho^2\mathbf{E}\left[\left\|\mathbf{x}^k\right\|_W^2\right]-\rho\mathbf{E}\left[\left\|\mathbf{x}^{k+1}-\mathbf{x}^k\right\|_{\mathcal{A}_{\beta,\eta}+\rho W-R}^2\right] \nonumber\\
		&-\rho\left((1-\eta)\gamma+\eta+\beta\right)\mathbf{E}\left[\left\|\mathbf{x}^k-\mathbf{x}^*\right\|^2\right] \nonumber\\
		&-\rho(\eta+\frac{1-\eta}{\gamma})G_2^2-2(G_1^2+G_2^2). \label{eq_lemma2}
		\vspace{-0.1cm}
	\end{align}
	By the strong convexity of $f$, (\ref{eq_lemma3_2}) holds.
\end{proof}

Next, we provide an upper bound on $\mathbf{E}\left[\left\|\mathbf{z}^k-\mathbf{z}^*\right\|_Q^2\right]$. For any $c_1,c_2>0$, combining (\ref{zeroth_primal_backtracking_iteration}), (\ref{dual_optimum}), (\ref{eq_lemma2_2}), the smoothness of $f$ and the fact that $\mathbf{q}^k=W^{\frac{1}{2}}\mathbf{v}^k$ together, we have 
\begin{align*}
	&\mathbf{E}[\|\mathbf{v}^k-\mathbf{v}^*\|^2]=\mathbf{E}[\|(W^\dagger)^{\frac{1}{2}}W^{\frac{1}{2}}(\mathbf{v}^k-\mathbf{v}^*)\|^2]\\
	=&\mathbf{E}[\|(W^\dagger)^{\frac{1}{2}}((A^k)^{-1} H^k(\mathbf{x}^k-\mathbf{x}^{k+1})-\rho W\mathbf{x}^k-\widetilde{g}_{\mu}(\mathbf{x}^k) \notag\\
	&+\nabla f(\mathbf{x}^*))\|^2] \notag\\
	\leq&(1+c_1)\mathbf{E}[\|(W^\dagger)^{\frac{1}{2}}((A^k)^{-1} H^k(\mathbf{x}^k-\mathbf{x}^{k+1})-\widetilde{g}_\mu(\mathbf{x}^k) \notag\\
	&+\nabla f(\mathbf{x}^k))\|^2] \notag\\
	&+(1+\frac{1}{c_1})\mathbf{E}\left[\|(W^\dagger)^{\frac{1}{2}}(\rho W\mathbf{x}^k+\nabla f(\mathbf{x}^k)-\nabla f(\mathbf{x}^*))\|^2\right]\\
	\leq&\frac{2(1+c_1)}{\lambda_W}\mathbf{E}\left[\|\mathbf{x}^{k+1}-\mathbf{x}^k\|_{\left(H^k(A^k)^{-1}\right)^2}^2\right]\\
	&+\rho^2(1+\frac{1}{c_1})(1+c_2)\mathbf{E}\left[\|\mathbf{x}^k\|_W^2\right]+\frac{2(1+c_1)(G_1+G_2)^2}{\lambda_W}\\
	&+\frac{(1+1/c_1)(1+1/c_2)}{\lambda_W}\mathbf{E}\left[\|\mathbf{x}^{k}-\mathbf{x}^*\|_{\Lambda_M^2}^2\right].
\end{align*}
\vspace{-0.2cm}
Since $\left(H^k(A^k)^{-1}\right)^2\preceq\left(\underline{\alpha}^{-1}\right)^2(\frac{1}{\theta}\Lambda_M+D)^2$, we then have 
\begin{align}
	\vspace{-0.1cm}
	&\mathbf{E}\left[\left\|\mathbf{z}^k-\mathbf{z}^*\right\|_Q^2\right] \notag\\
	\leq&\frac{2(1+c_1)}{\lambda _W}\mathbf{E}\left[\left\|\mathbf{x}^{k+1}- \mathbf{x}^k\right\|_{\left(\underline{\alpha}^{-1}\right)^2(\frac{1}{\theta}\Lambda_M+D)^2}^2\right] \notag\\
	&+\rho^2(1+\frac{1}{c_1})(1+c_2)\mathbf{E}\left[\left\|\mathbf{x}^k\right\|_W^2\right] \notag\\
	&+\mathbf{E}\left[\left\|\mathbf{x}^k-\mathbf{x}^*\right\|^2_\mathcal{B}\right]+\frac{2(1+c_1)(G_1+G_2)^2}{\lambda_W}.\label{lemma4_1}
	\vspace{-0.1cm}
\end{align}

Then, note from (\ref{eq_lemma3_2}) and (\ref{lemma4_1}) that for any $\delta\in(0,1)$,
\vspace{-0.1cm}\begin{align}
	&(1-\delta)\mathbf{E}\left[\left\|\mathbf{z}^k-\mathbf{z}^*\right\|_Q^2\right]-\mathbf{E}\left[\left\|\mathbf{z}^{k+1}-\mathbf{z}^*\right\|_Q^2\right] \notag\\
	\geq&\rho\lambda_{\min}\left(\delta_c\mathbf{I}_{Nd}-\frac{\delta\mathcal{B}}{\rho}\right)\mathbf{E}\left[\left\|\mathbf{x}^k-\mathbf{x}^*\right\|^2\right] \notag\\
	&+(\rho^2(1-\delta(1+1/c_1)(1+c_2))\mathbf{E}\left[\left\|\mathbf{x}^k\right\|_W^2\right]-G \notag\\
	&-\mathbf{E}\left[\left\|\mathbf{x}^{k+1}-\mathbf{x}^k\right\|_{\frac{2\underline{\alpha}^{-2}\delta(1+c_1)}{\lambda_W}(\frac{1}{\theta}\Lambda_M+D)^2+\rho(\mathcal{A}_{\beta,\eta}+\rho W-R)}^2\right]. \nonumber
\end{align}
\vspace{-0.02cm}
In order to guarantee $\delta_c>0$ for some $\beta$, $\eta$ and $\gamma$, we need to bound these parameters more strictly. To let $2m(1-\eta)-\eta-\beta>0$, we impose $\eta>1$ and $\beta>\underline{\alpha}^{-1}>0$. Then, set $\gamma>\frac{2m(\eta-1)+\eta+\beta}{\eta-1}$ such that $(\gamma-2m)(\eta-1)>\eta+\beta$. To satisfy (\ref{theorem1_1}), it suffices to let
\vspace{-0.1cm}
\begin{align}
	\vspace{-0.15cm}
	&\delta_c\mathbf{I}_{Nd}-\frac{\delta\mathcal{B}}{\rho}\succeq \mathbf{O}_{Nd},\label{theorem1_4}\\
	&\rho^2(1-\delta(1+c_2)(1+1/c_1))\geq0,\label{theorem1_5}\\
	&\frac{2\underline{\alpha}^{-2}\delta(1+c_1)}{\lambda_W}(\frac{1}{\theta}\Lambda_M+D)^2+\rho(\mathcal{A}_{\beta,\eta}+\rho W-R)\preceq\mathbf{O}_{Nd} \label{theorem1_6}.
	\vspace{-0.1cm}
\end{align}

From (\ref{theorem1_4})-(\ref{theorem1_6}), to guarantee the existence of $\delta\in(0,1)$, we need $\kappa_{\beta,\eta}>0$. By (\ref{eq_lemma3_1}) and the fact that $\beta>\underline{\alpha}^{-1}$ and $\underline{A}\preceq\mathbf{I}_{Nd}$, we have $\mathcal{A}_{\beta,\eta}+\rho W \prec R$ and thus $\kappa_{\beta,\eta}>0$ always exists. In conclusion, $\delta\in(0,1)$ satisfying (\ref{theorem1_4})-(\ref{theorem1_6}) is given by (\ref{theorem1_2}). Finally, from (\ref{theorem1_1}), for each $k\geq0$ we have
\vspace{-0.2cm}
\begin{align*}
	\mathbf{E}\left[\left\|\mathbf{z}^k-\mathbf{z}^*\right\|_Q^2\right]&\leq(1-\delta)^k\mathbf{E}\left[\left\|\mathbf{z}^0-\mathbf{z}^*\right\|_Q^2\right]+\sum_{t=0}^{k-1}(1-\delta)^tG.
	\vspace{-0.22cm}
\end{align*}
By taking the limit $k\rightarrow\infty$, we show (\ref{theorem1_7}) holds.


\bibliographystyle{IEEEtran}

\end{document}